\documentclass[twoside,11pt,reqno]{amsart}

\usepackage{upref,amsxtra,amssymb,amscd}
\usepackage{varioref}
\usepackage{verbatim}
\usepackage{epsfig}
\usepackage{color}

\usepackage{epic,eepic,eucal}
\usepackage{enumerate}

\def\Ba{        \textrm{\bf Bad}}

\def\t{           \theta}

\def\a{         \alpha}
\def\b{         \beta}

\newcommand{\NN}{{\mathbb N}}
\newcommand{\RR}{{\mathbb R}}
\newcommand{\TT}{{\mathbb T}}

\newcommand{\ZZ}{{\mathbb Z}}

\newtheorem{theo}{\sc Theorem}[section]

\newtheorem{lemm}[theo]{\sc Lemma}
\newtheorem{coro}[theo]{\sc Corollary}

\theoremstyle{definition}

\theoremstyle{remark}

\newtheorem{rema}[theo]{\sc Remark}

\numberwithin{equation}{section}

\begin{document}
\title{Badly approximable affine forms and Schmidt games}
\author{Jimmy Tseng}
\address{Jimmy Tseng, Department of Mathematics, The Ohio State University, Columbus, OH 43210}
\email{tseng@math.ohio-state.edu}

\begin{abstract}
For any real number $\t$, the set of all real numbers $x$ for which there exists a constant $c(x) > 0$ such that $\inf_{p \in \ZZ} |\t q - x - p| \geq \frac{c(x)}{|q|}$ for all $q \in \ZZ \backslash \{0\}$ is an $1/8$-winning set.
\end{abstract}

%\date{\today}
\maketitle
\section{Introduction}\label{secIntro}
Let $M_{m,n}(\RR)$ denote  the set of $m \times n$ real matrices and $\widetilde{M}_{m,n}(\RR)$ denote $M_{m,n}(\RR) \times \RR^m$.  The element in $\widetilde{M}_{m,n}(\RR)$ corresponding to $A \in M_{m,n}(\RR)$ and $\mathbf{b} \in \RR^m$ will be expressed as $\langle A, \mathbf{b} \rangle$.  Consider the following well-known sets from the theory of Diophantine approximation~\cite{Kl}:  \[\Ba(m,n) := \Big\{\langle A, \mathbf{b} \rangle \in \widetilde{M}_{m,n}(\RR)  \mid \exists \ c(A, \mathbf{b}) > 0 \textrm{ s.t. } \|A \mathbf{q} -\mathbf{b} \|_{\ZZ} \geq \frac{c(A, \mathbf{b})}{\|\mathbf{q}\|^{n/m}} \ \forall \mathbf{q} \in \ZZ^n \backslash \{\mathbf{0}\}\Big\}\] where $\|\cdot\|$ is the sup norm on $\RR^k$ and $\| \cdot \|_{\ZZ}$ is the norm on $\RR^k$ given by $\|\mathbf{x}\|_{\ZZ} := \inf_{p \in \ZZ^k}\|\mathbf{x} - \mathbf{p}\|$.  The set $\Ba(m,n)$ is called the \textbf{set of badly approximable systems of} $\mathbf{m}$\textbf{ affine forms in} $\mathbf{n}$ \textbf{variables.}  For any $\mathbf{b} \in \RR^m$, let $\Ba^{\mathbf{b}}(m,n) := \{A \in M_{m,n}(\RR) \mid \langle A, \mathbf{b} \rangle \in \Ba(m,n)\}$, and, for any $A \in M_{m,n}(\RR)$, let $\Ba_A(m,n):= \{\mathbf{b} \in \RR^m \mid  \langle A, \mathbf{b} \rangle \in \Ba(m,n)\}$.

The set $\Ba^{\mathbf{0}}(m,n)$ is called the \textbf{set of badly approximable systems of} $\mathbf{m}$\textbf{ linear forms in} $\mathbf{n}$ \textbf{variables} and is an important and classical object of study in the theory of Diophantine approximation.  Although it is a Lebesgue null set (Khintchine, 1926), it has full Hausdorff dimension and, even stronger, is winning (Schmidt, 1969).  Winning sets have a few other properties besides having full Hausdorff dimension; see Subsection~\ref{subsecBWSCF} for more details. 

For the larger set $\Ba(m,n)$, however, less is known.  Among its known properties are that it has Lebesgue measure zero, but full Hausdorff dimension.  The former property follows from the doubly metric inhomogeneous Khintchine-Groshev Theorem (\cite{Ca}, Chapter VII, Theorem II).  The latter property is a result of D. Kleinbock (1999) proved using mixing of flows on the space of lattices~\cite{Kl}.  Recently (2008), Y.~Bugeaud, S.~Harrap, S.~Kristensen, and S.~Velani have given a simpler proof of Kleinbock's result; their main result is that, for every $A$, $\Ba_A(m,n)$ (and some related sets) has full Hausdorff dimension~\cite{BHKV}.  Using the Marstrand slicing theorem (\cite{Fal}, Theorem 5.8), Kleinbock's result follows.  In view of these results, a natural question that arises is whether, like $\Ba^{\mathbf{0}}(m,n)$, these sets $\Ba_A(m,n)$ and $\Ba(m,n)$ are winning instead of just having full Hausdorff dimension.  In this note, we show that $\Ba_\t(1,1)$ is winning for every real number $\t$.\footnote{For $\Ba_\t(1,1)$, we have a slight strengthening of the aforementioned consequence of the Khintchine-Groshev Theorem:  $\Ba_\t(1,1)$ has Lebesgue measure zero for every irrational number $\t$~\cite{Ki}.  This result is essentially a corollary of two elementary facts from the theory of continued fractions (see~\cite{T1} for this short, second proof and for a connection with shrinking targets).  There is yet a third proof of this result; see~\cite{BBDV}.}  For results and open questions concerning general $n$ and $m$, see Remark~\ref{secConcl} below.

\subsection{Statement of results}  Our main result, which generalizes the $m=n=1$ case of the aforementioned main result in~\cite{BHKV} (their main result is Theorem 1 of~\cite{BHKV}), is the following (see Subsection~\ref{subsecBWSCF} for the definition of $1/8$-winning):

\begin{theo} \label{thmDualBAwinning}
For any real number $\t$, $\Ba_\t(1,1)$ is an $1/8$-winning set.
\end{theo}

\noindent This theorem is proved in Section~\ref{secPTDualBAwinning} below.  A number of corollaries will follow immediately because of properties of winning sets (see Subsection~\ref{subsecBWSCF}).  A model one is

\begin{coro}\label{coroAffineImageWinning}
For any countable set $\{\t_n\} \subset \RR$ and any countable family $\{f_m\}$ of invertible affine maps $\RR \rightarrow \RR$, the set $\cap_{m=1}^\infty\cap_{n=1}^\infty f_m(\Ba_{\t_n}(1,1))$ is $1/8$-winning and thus has full Hausdorff dimension.\end{coro}

\subsection{Background on winning sets and continued fractions}\label{subsecBWSCF}

The proof of our result requires two tools:  Schmidt games (see~\cite{Sch2} for a reference) and continued fractions (see~\cite{Kh} for a reference).  We will discuss both.  

W. Schmidt introduced the games which now bear his name in~\cite{Sch2}.  Let $0 < \a <1$ and $0 < \b<1$.  Let $S$ be a subset of a complete metric space $M$.  Two players, Black and White, alternate choosing nested closed balls $B_1 \supset W_1 \supset B_2 \supset W_2 \cdots$ on $M$.  The radius of $W_n$ must be $\a$ times the radius of $B_n$, and the radius of $B_n$ must be $\b$ times the radius of $W_{n-1}$.  The second player, White, wins if the intersection of these balls lies in $S$.  A set $S$ is called \textbf{$(\a, \b)$-winning} if White can always win for the given $\a$ and $\b$.  A set $S$ is called {\bf $\a$-winning} if White can always win for the given $\a$ and any $\b$.  A set $S$ is called {\bf winning} if it is $\a$-winning for some $\a$.  Schmidt games have four important properties for us~\cite{Sch2}: \medskip

\noindent$\bullet$ The sets in $\RR^n$ which are $\a$-winning have full Hausdorff dimension.

\noindent$\bullet$ Countable intersections of $\a$-winning sets are again $\a$-winning.

\noindent$\bullet$ The bilipschitz image of an $\a$-winning set is $\a$-winning.

%\noindent$\bullet$ If a set is $\a$-winning, then it is also $\a'$-winning for all $0 < \a' \leq \a$.
\noindent$\bullet$ Let $0 < \a \leq 1/2$.  If a set in a Banach space of positive dimension is $\a$-winning, then the set with a countable number of points removed is also $\a$-winning. \medskip

Let us now discuss continued fractions.  Let $p_i/q_i$ be the $i$-th order convergent of an irrational number $\t$.   Define $$\Delta_i := \|\t q_i\|_{\ZZ}.$$  We will use the following well-known facts: \medskip

%\noindent$\bullet$ For $i \geq 2$, $q_i \geq 2 q_{i-2}$.

\noindent$\bullet$ For all $i \in \NN$, $\frac 1 2 \Delta_{i-1}^{-1} < q_i \leq \Delta_{i-1}^{-1}$. 

\noindent$\bullet$ Let $0 \leq j < k <q_i$.  Then, $\|\t k-\t j\|_{\ZZ} \geq \Delta_{i-1}$.

\subsection{The setup}
Let $\t \in \RR$.  Define \[\Ba^+_\t := \Big\{x \in \RR \mid \exists \ c(x) > 0 \textrm{ s.t. } \|\t q - x\|_{\ZZ} \geq \frac{c(x)}{q} \ \forall q \in \NN\Big\}.\]  Note that $\Ba_\t(1,1) = \Ba^+_\t \cap -\Ba^+_\t$; thus showing $\Ba^+_\t$ is $1/8$-winning will prove Theorem~\ref{thmDualBAwinning}.  Also, we may assume that these sets are restricted to the circle $\TT^1:=\RR / \ZZ$, as they are invariant under integral translations.

Henceforth, let us consider $\Ba_\t^+$.  If $\t$ is rational, then the set is just $\TT^1$ with a finite number of points removed and hence is winning.  Therefore, we assume that $\t$ is irrational henceforth.

For convenience, let us call the elements in \[\{\t q  \in \TT^1 \mid q_{i} \leq q < q_{i+1}\}\] the \textbf{elements of generation $i$.}

\noindent Finally, we note a simple property of continued fractions:

\begin{lemm} \label{lemmCircleMultLaw}
Let $q_{i+1} \leq q < q_{i+2}$.  Given a $0< r < 1/2$ such that, for all elements $\t p$ of generations $\leq i$, $\|\t q - \t p\|_\ZZ \geq r \Delta_i$, then $ q \geq \frac r 2 q_{i+2}$.
\end{lemm}

\begin{proof}
There are unique numbers $0 \leq s < q_{i+1}$ and $1 \leq n \leq \lfloor \frac {q_{i+2}} {q_{i+1}} \rfloor$ such that $q = n q_{i+1} + s$.  Thus, $n \Delta_{i+1} = \|\t q - \t s \|_{\ZZ} \geq r \Delta_i$.  Hence, $q \geq r \frac{\Delta_i}{\Delta_{i+1}} q_{i+1} \geq \frac r 2 q_{i+2}.$
\end{proof} 

\section{A proof of Theorem~\ref{thmDualBAwinning}} \label{secPTDualBAwinning}

Let $\a = 1/8$ and $c = (\frac{(\a \b)}{4})^3.$  We will play an $(\a,\b)$-game on $\TT^1$.  Let us start with the following lemma, which tells us how to choose $W_m$ given $B_m$ (note that the radius of a ball $B$ will be denoted $\rho(B)$):  

\begin{lemm} \label{lemmGameStep}
Let $U$ be any union of balls on $\TT^1$ with radius $\leq (\a \b) \Delta_N/4$ around the elements of generations $\leq N$.  If \[(\a \b) \Delta_N< 2 \rho(B_m) \leq \Delta_N,\] then one can choose $W_m$ disjoint from $U$.
\end{lemm}

\begin{proof}
\medskip\noindent\textbf{Case:  $B_m$ does not intersect any ball of $U$.}  

Pick any allowed $W_m$.

\medskip\noindent\textbf{Case:  $B_m$ intersects exactly one ball of $U$.}  

Even if $B_m$ contains the whole ball of $U$, there is, at least, a subinterval in $B_m$ of length $1/4$ of the length of $B_m$ that misses $U$.  Pick $W_m$ to be in this subinterval.

\medskip\noindent\textbf{Case:  $B_m$ intersects more than one ball of $U$.}   Note that $B_m$ cannot intersect more than one element of generations $\leq N$ (unless one has exactly two elements of generations $\leq N$, one at each end).  Thus, at least a subinterval in $B_m$ of length $(1 - (\a \b)/2 ) \Delta_N \geq 1/2 \Delta_N$ does not meet $U$. Now $\a 2 \rho(B_m) \leq 1/8 \Delta_N$.  Therefore, we can choose $W_m$ to be in this subinterval.\end{proof}

Since the Schmidt game can be played until, for some $J \in \NN$, $2 \rho(B_J) \leq \Delta_1$, we may assume without loss of generality that $J = 1$.  Note that there exists a $N \geq 2$ such that $2 \rho(B_1) \leq \Delta_{N-1}$, but that $2 \rho(B_1) > \Delta_{N}$ (follows since $\Delta_{N} < \Delta_{N-1}$).

Also, there exists a $n_0 \in \NN$ such that $2 (\a \b)^{n_0-1} \rho(B_1) > \Delta_N$ and $2 (\a \b)^{n_0}   \rho(B_1) \leq \Delta_N$.  Thus, \begin{equation}\label{eqnbestfit}
(\a \b) \Delta_N < 2 (\a \b)^{n_0} \rho(B_1) \leq \Delta_N.
\end{equation}   We require that $N$ be the largest natural number for which (\ref{eqnbestfit}) holds.   

We intend to use induction.  In the initial induction step, consider the disjoint union of balls around each element of generations $\leq N$ of radius $(\a \b) \Delta_N/4$; call this union $U$.  By Lemma~\ref{lemmGameStep}, we may pick $W_{n_0+1}$ to miss $U$.  For any other step of the induction, $W_{n_0+1}$ is already chosen.

\subsection{Case: $\a \b \Delta_N > \Delta_{N+1}$.}
The condition implies that there exists a $n_1 \in \NN$ such that \[ (\a \b) \Delta_{N+1} < 2(\a \b)^{n_0+n_1}  \rho(B_1)\leq \Delta_{N+1}.\]  Also, there exists a  maximal $M \geq 1$ such that \[ (\a \b) \Delta_{N+M} < 2(\a \b)^{n_0+n_1}  \rho(B_1)\leq \Delta_{N+M}.\]  Moreover, $(\a \b) \Delta_{N+1} < \Delta_{N+M}.$

For any element $\t q$ of generation $N+1$ in $W_{n_0+1}$, $q \geq \frac{(\a \b)} {8} q_{N+2}$ by Lemma~\ref{lemmCircleMultLaw}.  For any element $\t q$ of generations $> N+1$ in $W_{n_0+1}$, it is obvious that $q \geq \frac{(\a \b)} {8} q_{N+2}$.
%$c = (\a \b) (1 - \frac{\a \b} {8}) / 4$.
Thus, for all such $\t q$, \[\frac c q \leq \frac {(\a \b)^2 \Delta_{N+1}} 4 \leq \frac {(\a \b) \Delta_{N+M}} 4.\]  Now play freely until $B_{n_0+n_1+1}$ is chosen.  Again by Lemma~\ref{lemmGameStep}, we can choose $W_{n_0+n_1+1}$ to miss the balls of radius $ (\a \b) \Delta_{N+M}/ 4$ around the elements of generations $N+1$ to $N+M$.  

\subsection{Case:  $\a \b \Delta_N \leq \Delta_{N+1}$}  It is easy to see from the theory of continued fractions that there exist a $K \in \NN$ such that $(\a \b) \Delta_n > \Delta_{n+K}$ for all $n \in \NN$.  Therefore, the condition implies that there exists a $1 \leq m \leq K-1$ such that \[\Delta_{N+m+1} < \a \b \Delta_N \leq \Delta_{N+m}.\]  Thus, we have \[(\a \b)^2 \Delta_{N+m} < (\a \b)^2 \Delta_N < 2 (\a \b)^{n_0+1} \rho(B_1)  \leq \a \b \Delta_N \leq \Delta_{N+m}.\]  %\textrm{ and } 2 (\a \b)^{n_0+1} \rho(B_1)
If $(\a \b)^2 \Delta_{N+m} <  2 (\a \b)^{n_0+1} \rho(B_1) \leq (\a \b) \Delta_{N+m},$ then \[(\a \b) \Delta_{N+m} <  2 (\a \b)^{n_0} \rho(B_1) \leq \Delta_{N+m}.\]  Since $N$ is the largest natural number for which (\ref{eqnbestfit}) holds, we obtain that $m=0$, a contradiction.  %Let us, for convenience, denote $n:=n_0$ for this possibility.  

Thus, we must conclude that \[(\a \b) \Delta_{N+m} <  2 (\a \b)^{n_0+1} \rho(B_1) \leq \Delta_{N+m}.\]  Now, there exists a $n_1 \in \NN$ such that \[(\a \b) \Delta_{N+m+1} < 2 (\a\b)^{n_0+n_1} \rho(B_1) \leq \Delta_{N+m+1}.\]   Also, there exists a maximal $M \in \NN$ such that \[(\a \b) \Delta_{N+m+M} < 2 (\a\b)^{n_0+n_1} \rho(B_1) \leq \Delta_{N+m+M}.\]  Moreover, $(\a \b) \Delta_{N+m+1} < \Delta_{N+m+M}$. 

If $n_1=1$, then even more is true:  $(\a \b) \Delta_{N+m} < \Delta_{N+m+M}$.  Now note that, for the elements $\t q$ of generations $N+1$ to $N+m+M$, we have \[\frac c q \leq \frac {c} {q_{N+1}} \leq \frac{(\a \b) \Delta_{N+m+M}}4.\]  Consider the disjoint union of balls around each element of generations $\leq N+m+M$ of radius $(\a \b) \Delta_{N+m+M}/4$; call this union $U$. Again by Lemma~\ref{lemmGameStep}, we can pick $W_{n_0+2}$ to miss $U$.

Otherwise, $n_1 \geq 2$.  Now note that, for the elements $\t q$ of generations $N+1$ to $N+m$, we have \[\frac c q \leq \frac {c} {q_{N+1}} \leq \frac{(\a \b) \Delta_{N+m}}4.\] Consider the disjoint union of balls around each element of generations $\leq N+m$ of radius $(\a \b) \Delta_{N+m}/4$; call this union $U$. Again by Lemma~\ref{lemmGameStep}, we can pick $W_{n_0+2}$ to miss $U$.

For any element $\t q$ of generation $N+m+1$ in $W_{n_0+2}$, $q \geq \frac{(\a \b)} {8} q_{N+m+2}$ by Lemma~\ref{lemmCircleMultLaw}.  
%$c = (\a \b)^2 (1 - \frac{(\a \b)^2} {8}) / 4$.  
For any element $\t q$ of generations $> N+m+1$ in $W_{n_0+2}$, it is obvious that $q \geq \frac{(\a \b)} {8} q_{N+m+2}$.  
Thus, for all such $\t q$, \[\frac c q \leq \frac {(\a \b)^2 \Delta_{N+m+1}} 4 \leq \frac {(\a \b) \Delta_{N+m+M}} 4.\]  Now play freely until $B_{n_0+n_1+1}$ is chosen.  Again by Lemma~\ref{lemmGameStep}, we can choose $W_{n_0+n_1+1}$ to miss the balls of radius $(\a \b) \Delta_{N+m+M} / 4$ around the elements of generations $N+m+1$ to $N+m+M$.  

\bigskip

Using these two cases inductively, one can show that the set \[\Big\{x \in \RR \mid \exists \ c(x) > 0 \textrm{ s.t. } \|\t q - x\|_{\ZZ} \geq \frac{c(x)}{q} \ \forall q \geq q_{N+1}\Big\}\] is  $1/8$-winning.  By shrinking $c(x)$ for each $x$, we note that this set is $\Ba^+_\t$.  The proof is complete.

\begin{rema}
If $\t$ is a badly approximable number\footnote{In our notation, $\t  \in \Ba^0(1,1)$.}, one can easily see from the continued fraction expansion of $\t$ that there exists an upper bound for $\Delta_n / \Delta_{n+1}$ independent of $n$.  This uniform bound allows us to simplify the above proof for $\t$ badly approximable (however, we conclude that the set is $\a$-winning for an $\a$ depending on this uniform bound).
\end{rema}

\begin{rema}\label{secConcl}
It should be noted that Manfred Einsiedler and the author, in joint work in preparation, are able to generalize Theorem 1 of~\cite{BHKV}, using a method different from the one presented in this note, to conclude winning, instead of just having full Hausdorff dimension.  Thus, as a special case, we can show that $\Ba_A(m,n)$ is winning for every $A \in M_{m,n}(\RR)$.  Whether $\Ba(m,n)$ is winning, however, is still an open question.  See Section 5 of~\cite{Kl} for a list of this and other open questions related to winning.\end{rema}

\section*{Acknowledgements} 

The author would like to thank Manfred Einsiedler and Dmitry Kleinbock for their helpful comments.  The author is particularly grateful to Dmitry Kleinbock for pointing out an improvement to the statement of Corollary~\ref{coroAffineImageWinning}.


\begin{thebibliography}{99}

\bibitem{BBDV} V. Beresnevich, V. Bernik, M. Dodson, and S. Velani, {\em Classical metric Diophantine approximation revisited,} preprint, arXiv:0803.2351v1 (2008).

\bibitem{BHKV} Y. Bugeaud, S. Harrap, S. Kristensen, and S. Velani, {\em On shrinking targets for $\ZZ^m$ actions on tori,} preprint, arXiv:0807.3863v1 (2008).

\bibitem{Ca} J. Cassels, ``An Introduction to Diophantine Approximation,'' Cambridge Tracts in Mathematics and Mathematical Physics {\bf 45}, Cambridge
University Press, Cambridge, UK, 1957.

\bibitem{Fal} K. Falconer, ``The geometry of fractal sets," Cambridge Tracts in Mathematics, {\bf 85}, Cambridge University Press, Cambridge, 1986.

\bibitem{Kh} A. Khinchin, ``Continued Fractions,'' The University of Chicago Press, Chicago, 1964.

\bibitem{Ki} D. Kim, {\em The shrinking target property of irrational rotations}, Nonlinearity {\bf 20} (2007), 1637--1643.

\bibitem{Kl} D. Kleinbock, {\em Badly approximable systems of affine forms}, J. Number Theory {\bf 79} (1999), 83 --102. 

\bibitem{Sch2} W. Schmidt, {\em Badly approximable numbers and certain games,} Trans. Amer. Math. Soc. {\bf 123} (1966), 178--199.

\bibitem{T1} J. Tseng, {\em On circle rotations and the shrinking target properties,} Discrete Contin. Dyn. Syst. {\bf 20} (2008), 1111-1122.

\end{thebibliography}
\end{document}